\documentclass[12pt]{article}

\usepackage{amsmath,amssymb,amsfonts,amsthm,epic}

\usepackage{graphicx}

\theoremstyle{plain}

\newtheorem{ex}[subsubsection]{Example}

\theoremstyle{definition}
\newtheorem{rem}[subsubsection]{Remark}
\newtheorem{defn}[subsubsection]{Definition}
\newtheorem{nota}[subsubsection]{Notation}

\newcommand\brem{\begin{rem}\begin{sffamily}\begin{upshape}}
\newcommand\erem{\end{upshape}\end{sffamily}\end{rem}}
\newcommand\bdefn{\begin{defn}\begin{rm}}
\newcommand\edefn{\end{rm}\hfill$\Box$\end{defn}}
\newcommand\bnot{\begin{nota}\begin{rm}}
\newcommand\enot{\end{rm}\hfill$\Box$\end{nota}}
\newcommand\bex{\begin{ex}\begin{rm}}
\newcommand\eex{\end{rm}\hfill$\Box$\end{ex}}

\numberwithin{equation}{subsection}
\title{A rank computation problem related to the HK function of trinomial hypersurfaces}
\date{}
\author {Shyamashree Upadhyay\\Department of Mathematics\\ Indian Institute of Technology, Guwahati\\Assam-781039, INDIA\\email:shyamashree@iitg.ernet.in}
\begin{document}
\maketitle

\begin{abstract}
In this article, I provide a solution to a rank computation problem related to the computation of the Hilbert-Kunz function for any disjoint-term trinomial hypersurface, over any field of characteristic $2$. This rank computation problem was posed in \cite{Su1}. The formula for the Hilbert-Kunz function for disjoint-term trinomial hypersurfaces can be deduced from the result(s) of this article, modulo a lot of tedious notation.
\end{abstract}
\tableofcontents
\section{Introduction}\label{s.Introduction}
A `disjoint-term trinomial hypersurface' is defined in \S 2 of \cite{Su1}. In the last section of \cite{Su1}, I gave an algorithm for computing the Hilbert-Kunz function for any disjoint-term trinomial hypersurface in general, over any field of arbitrary positive characteristic. But the algorithm was given modulo a rank computation problem for $3$ types of systems of linear equations which were given in tables $5,6$ and $8$ of \cite{Su1}. In this article, I provide a solution to this rank computation problem over fields of characteristic $2$. It will be nice if one can provide a similar solution over fields of arbitrary positive characteristic $p$.

The formula for the Hilbert-Kunz function for disjoint-term trinomial hypersurfaces follows as a corollary to the solution of this rank computation problem, but a lot of tedious notation is needed for deducing the formula. Hence in this article, I am avoiding the work on the deduction of the Hilbert-Kunz function from this. 

In the article \cite{Su1}, I was suspecting that due to the weird combinatorial pattern of the matrices in tables $5,6$ and $8$, there can be examples of disjoint-term trinomial hypersurfaces for which the corresponding Hilbert-Kunz multiplicity can be irrational. But my suspicion turns out to be incorrect (see the papers \cite{Mo-Teix2} and \cite{Mo-Teix1} for a reasoning). However if we consider trinomial hypersurfaces which are defined by a polynomial not having disjoint terms in it, then we will encounter a similar-looking (but larger and more complicated) rank computation problem for those. And I suspect that the Hilbert-Kunz multiplicity may become irrational for those hypersurfaces because of the more complicated nature of the system of linear equations involved in it. In fact, the solution to the rank computation problem given in this article is also important for the next article in which I will consider trinomial hypersurfaces which are defined by a polynomial not having disjoint terms in it, because something similar happens there. 

\section{Stating the problem in a general way}\label{s.statetheproblem}
The following combinatorial fact is well known (called the \textit{Chu-Vandermonde identity}):
\begin{center}
For any two natural numbers $m$ and $n$ and for any non-negative integer $j$, the standard inner product of the two vectors $(\ ^{m}C_0,\ ^{m}C_1,\ldots,\ ^{m}C_j)$ and $(\ ^{n}C_j,\ ^{n}C_{j-1},\ldots,\ ^{n}C_0)$ equals $\ ^{m+n}C_j$ where we declare that for any positive integer $M$, $\ ^{M}C_l=0$ if $l>M$.
\end{center}

It follows from the above fact and from the pattern of the matrices present in tables $5,6$ and $8$ of the article \cite{Su1} that
\begin{itemize}
\item for the linear systems corresponding to the matrices in tables $5$ and $8$, we need to solve the general combinatorial problem mentioned below as \textbf{Problem I} and\\
\item for the linear systems corresponding to the matrices in tables $6$, we need to solve the general combinatorial problem mentioned below as \textbf{Problem II}. 
\end{itemize}
It is worthwhile to mention here that in this article, we will be considering systems of linear equations over fields of characteristic $2$ only.
\subsection{Statement of Problem I}\label{ss.state-ProblemI}
Given any positive integer $M$ and any non-negative integer $j$, what all values of positive integers $k$ and $l$ solve the following system of $(l+1)$ many linear equations in $(k+1)$ many unknowns given by $AX=b$ where the coefficient matrix $A$ is the $(l+1)\times(k+1)$ matrix given below in table $1$ and $b$ is the $(l+1)\times(1)$ column vector $(0,0,\ldots,1)^{t}$. 
\begin{table}[hbt]
\begin{center}
\caption{The coefficient matrix $A$ of Problem I} 
\begin{tabular}{|c|c|c|c|c|c|c|}
\hline
$\ ^{M}C_{j+l+k}$ & $\ ^{M}C_{j+l+(k-1)}$ & $\cdots$ & $\cdots$ & $\cdots$ & $\ ^{M}C_{j+l+1}$ & $\ ^{M}C_{j+l}$ \\
\hline
$\ ^{M}C_{j+l+(k-1)}$ & $\diagup$ & $\diagup$ & $\diagup$ & $\diagup$ & $\diagup$ & $\ ^{M}C_{j+(l-1)}$ \\
\hline
$\diagup$ & $\diagup$ & $\diagup$ & $\diagup$ & $\diagup$ & $\diagup$ & $\vdots$ \\
\hline
$\diagup$ & $\diagup$ & $\diagup$ & $\diagup$ & $\diagup$ & $\diagup$ & $\ ^{M}C_{j+1}$ \\
\hline
$\diagup$ & $\diagup$ & $\diagup$ & $\diagup$ & $\diagup$ & $\diagup$ & $\ ^{M}C_j$ \\
\hline
\end{tabular}
\end{center}
\end{table}
The lines along the south-west$\leftrightarrow$north-east direction (in table $1$) indicate the continuation of the same entry along that direction.
\subsection{Statement of Problem II}\label{ss.state-ProblemII}
Given any positive integer $M$ and non-negative integer $j$, what all values of positive integers $k,l$ and $q$ solve the following system of $(l+q+1)$ many linear equations in $(k+1)$ many unknowns given by $AX=b$ where the coefficient matrix $A$ is the $(l+q+1)\times(k+1)$ matrix given below in table $2$ and $b$ is the $(l+q+1)\times(1)$ column vector $(0,0,\ldots,1)^{t}$. 
\begin{table}[hbt]
\begin{center}
\caption{The coefficient matrix $A$ of Problem II} 
\begin{tabular}{|c|c|c|c|c|c|c|}
\hline
$\heartsuit_q^0$ & $\heartsuit_q^1$ & $\cdots$ & $\cdots$ & $\cdots$ & $\heartsuit_q^{k-1}$ & $\heartsuit_q^{k}$ \\
\hline
$\vdots$ & $\vdots$ & $\vdots$ & $\vdots$ & $\vdots$ & $\vdots$ & $\vdots$ \\
\hline
$\heartsuit_2^0$ & $\heartsuit_2^1$ & $\cdots$ & $\cdots$ & $\cdots$ & $\heartsuit_2^{k-1}$ & $\heartsuit_2^{k}$ \\
\hline
$\heartsuit_1^0$ & $\heartsuit_1^1$ & $\cdots$ & $\cdots$ & $\cdots$ & $\heartsuit_1^{k-1}$ & $\heartsuit_1^{k}$ \\
\hline
$\ ^{M}C_{j+l+k}$ & $\ ^{M}C_{j+l+(k-1)}$ & $\cdots$ & $\cdots$ & $\cdots$ & $\ ^{M}C_{j+l+1}$ & $\ ^{M}C_{j+l}$ \\
\hline
$\ ^{M}C_{j+l+(k-1)}$ & $\diagup$ & $\diagup$ & $\diagup$ & $\diagup$ & $\diagup$ & $\ ^{M}C_{j+(l-1)}$ \\
\hline
$\diagup$ & $\diagup$ & $\diagup$ & $\diagup$ & $\diagup$ & $\diagup$ & $\vdots$ \\
\hline
$\diagup$ & $\diagup$ & $\diagup$ & $\diagup$ & $\diagup$ & $\diagup$ & $\ ^{M}C_{j+1}$ \\
\hline
$\diagup$ & $\diagup$ & $\diagup$ & $\diagup$ & $\diagup$ & $\diagup$ & $\ ^{M}C_j$ \\
\hline
\end{tabular}
\end{center}
\end{table}
In table $2$,
\begin{center}
$M=\alpha+\delta$ for some natural numbers $\alpha$ and $\delta$ (given to us).\\
For any $1\leq i\leq q$ and $0\leq r\leq k$, $\heartsuit_{i}^{r}:=$ the standard inner product\\ of the two vectors
$(\ ^{\alpha}C_0,\ ^{\alpha}C_1,\ldots,\ ^{\alpha}C_{j+k+l-r})$\\ and $(\ ^{\delta}C_{j+k+l+i-r},\ldots,\ ^{\delta}C_{i+1},\ ^{\delta}C_i)$.\\
The lines along the south-west$\leftrightarrow$north-east direction indicate the continuation of the same entry along that direction.
\end{center}
\section{Structure of the string of Binomial coefficients}\label{s.structure-binomial-string}
For providing solutions to both problems I and II, we first need to have an account of the position of all odd entries in the string of binomial coefficients of any given positive integer $M$ starting from $\ ^{M}C_0$ to $\ ^{M}C_M$. This is given as follows:
\begin{center}
Let $M=2^{l_m}+2^{l_{m-1}}+\cdots+2^{l_1}$ be the binary expansion of $M$ where $l_m>l_{m-1}>\cdots>l_1$ are non-negative integers.

Let $\mathfrak{D}_M:=$ the set of all possible $(m+1)$-tuples $(\alpha_0,\alpha_1,\ldots,\alpha_m)$ generated out of the set $\{0,2^{l_1},2^{l_2},\ldots,2^{l_m}\}$ which satisfy the following $3$ properties simultaneously:
\begin{itemize}
\item $\alpha_0=1$.
\item $\alpha_1\geq\alpha_2\geq\cdots\geq\alpha_m$.
\item If $\alpha_i\neq 0$ for some $i\in\{1,\ldots,m\}$, then $\alpha_i>\alpha_j$ for any $j\in\{i+1,\ldots,m\}$.
\end{itemize}
Clearly, the total number of elements in the set $\mathfrak{D}_M$ is $2^m$. Let us put an order $\precsim$ on $\mathfrak{D}_M$ by declaring that 
$$(\alpha_0,\alpha_1,\ldots,\alpha_m)\precsim(\alpha_0´,\alpha_1´,\ldots,\alpha_m´)$$ if and only if\\
$$\alpha_i\leq\alpha_i´\ \forall i\in\{0,1,\ldots,m\}.$$
Observe that $\precsim$ is a total order on the set $\mathfrak{D}_M$. The positions of the odd entries in the string of binomial coefficients of $M$ are indexed by elements of the set $\mathfrak{D}_M$; and are arranged (in ascending order) according to the order $\precsim$ on $\mathfrak{D}_M$. This indexing is given by:\\
For any $d=(\alpha_0,\alpha_1,\ldots,\alpha_m)\in\mathfrak{D}_M$, the position of the odd entry corresponding to it in the string of binomial coefficients of $M$ is $\Sigma_{i=0}^{m}\alpha_i$.\\

Since the total number of elements in the set $\mathfrak{D}_M$ is $2^m$, we can denote the elements of $\mathfrak{D}_M$ by $d_1^{(M)},d_2^{(M)},\ldots,d_{2^m}^{(M)}$ where $d_1^{(M)}\precsim d_2^{(M)}\precsim\ldots\precsim d_{2^m}^{(M)}$. And let us denote the corresponding positions of the odd entries in the string of binomial coefficients of $M$ by 
$$\Sigma d_1^{(M)},\Sigma d_2^{(M)},\ldots,\Sigma d_{2^m}^{(M)}.$$ 
\end{center}
Looking at the positions of all odd entries in the string of binomial coefficients of the given positive integer $M$, we can deduce some important facts about the position of all even entries in the string of binomial coefficients of $M$. These facts will be very much useful in writing down solutions to both the Problems I and II. These facts are as given below:

Let the binary expansion of $M$ be as mentioned above. Then 
\begin{enumerate}
\item the number of even entries between any two consecutive odd entries in the continuous string of binomial coefficients of $M$ is of the form $2^{l_{g}}-2^{l_{g-1}}-\ldots -2^{l_1}-1$ for some $g\in\{1,\ldots,m\}$.
\item If $l_g<l_{\hat{g}}$, then a bunch of $2^{l_{g}}-2^{l_{g-1}}-\ldots -2^{l_1}-1$ many even entries appears on both right and left sides of the bunch of $2^{l_{\hat{g}}}-2^{l_{\hat{g}-1}}-\ldots -2^{l_1}-1$ many even entries at equal distance from the central bunch.
\item Between any two consecutive bunches of $2^{l_{g}}-2^{l_{g-1}}-\ldots -2^{l_1}-1$ many even entries, there exists $(2^{g}-1)$ many bunches of even entries such that the number of even entries in each of these $(2^{g}-1)$ many bunches is different from $2^{l_{g}}-2^{l_{g-1}}-\ldots -2^{l_1}-1$.
\end{enumerate}
\section{Solutions to Problems I and II}\label{s.solutionprobIandII}
In subsection \ref{ss.solutionprobI} below, the notation and terminology will remain the same as in subsection \ref{ss.state-ProblemI}. Similarly, in subsection \ref{ss.solutionprobII} below, the notation and terminology will remain the same as in subsection \ref{ss.state-ProblemII}.
\subsection{Solution to Problem I}\label{ss.solutionprobI}
Problem I can be restated in the following way:---
\begin{center}
Given any positive integer $M$ and any non-negative integer $j$, for what all values of positive integers $k$ and $l$ does the string
$$\ ^{M}C_j, \ ^{M}C_{j+1},\ldots,\ldots,\ ^{M}C_{j+l+k}$$
of consecutive binomial coefficients of $M$ satisfy all the following $3$ properties simultaneously:
\begin{itemize}
\item $\ ^{M}C_j$ and $\ ^{M}C_{j+k+1}$ have different parity (that is, exactly one of them is even).
\item $\ ^{M}C_{j+t}=\ ^{M}C_{(j+t)+q(k+1)}$ for every $t\in\{1,\ldots,k+1\}$ and for every positive integer $q$ such that $(j+t)+q(k+1)\leq j+l+k$. In other words, the string of length $k+1$ beginning from $\ ^{M}C_{j+1}$ and ending at $\ ^{M}C_{j+k+1}$ should be repeated until we reach $\ ^{M}C_{j+l+k}$ (there can be a truncated string of length $k+1$ at the end).
\item  The string of length $k+1$ beginning from $\ ^{M}C_{j+1}$ and ending at $\ ^{M}C_{j+k+1}$ should contain \textit{evenly many} odd elements in it.
\end{itemize}  
\end{center}
There are many possible values of positive integers $k$ and $l$ which do the job. A list of all possible answers is given below in cases I and II.\\

\underline{\textbf{Case I:} When $\ ^{M}C_j$ is odd.}\\
Suppose that $\ ^{M}C_j$ equals $\Sigma d_i^{(M)}$ for some $i\in\{1,\ldots,2^m\}$, then we provide a solution to Problem I by considering different sub cases depending upon the value $\Sigma d_{i+1}^{M}-\Sigma d_{i}^{M}-1$ being $0$ or $>0$. Observe that $\Sigma d_{i+1}^{M}-\Sigma d_{i}^{M}-1$is nothing but the total number of zeroes lying in between the positions $\Sigma_{i+1}^{M}$ and $\Sigma_{i}^{M}$. \\

\underline{Sub case I.a:} When $\Sigma d_{i+1}^{(M)}-\Sigma d_{i}^{(M)}-1=0$.\\
If there exists an even natural number $s$ such that $\Sigma d_{i+s+1}^{(M)}-\Sigma d_{i+s}^{(M)}-1>0$, then take the string $\ ^{M}C_{j+1},\ldots,\ ^{M}C_{j+k+1}$ of length $(k+1)$ to be the continuous string of binomial coefficients of $M$ starting from position number $1+\Sigma d_i^{(M)}$ and ending at position number $\Sigma d_{i+s}^{(M)}+L$ (inclusive of the beginning and end points), where $L$ is some natural number such that $0<L\leq\Sigma d_{i+s+1}^{(M)}-\Sigma d_{i+s}^{(M)}-1$. Hence the natural number $k$ is given by $k=\Sigma d_{i+s}^{(M)}+L-\Sigma d_i^{(M)}-1$ where $L$ and $s$ are as mentioned above.

To determine the possible values of $l$ in this sub case, look at the entry in position number $\Sigma d_{i+s}^{(M)}+L+1$. If it is even, then take $l=1$. If it is odd, then look at the number of odd entries in the string of binomial coefficients of $M$ that come in continuum starting from position number $1+\Sigma d_i^{(M)}$ (including the odd entry at this position).  Let $p_j^{(M)}$ denote this number. Similarly, look at the number of odd entries in the string of binomial coefficients of $M$ that come in continuum starting from position number $1+L+\Sigma d_{i+s}^{(M)}$ (including the odd entry at this position).  Let $q_j^{(M)}$ denote this number. Tale $l$ to be any natural number such that $l\leq 1+min\{p_j^{(M)},q_j^{(M)}\}$.
\brem\label{r.subcaseI.a}
An even natural number $s$ as mentioned above in this sub case may not exist for some $M$ and for some position number $\Sigma d_i^{(M)}$. In that situation, this sub case is invalid.
\erem
\underline{Sub case I.b:} When $\Sigma d_{i+1}^{(M)}-\Sigma d_{i}^{(M)}-1>0$.\\
I will now list down the various possibilities under this sub case:
\begin{center}
Either\\
Take $k$ to be any natural number such that $k+1\leq\Sigma d_{i+1}^{(M)}-\Sigma d_{i}^{(M)}-1$. And take the string $\ ^{M}C_{j+1},\ldots,\ ^{M}C_{j+k+1}$ of length $(k+1)$ to be the continuous string of binomial coefficients of $M$ starting from position number $1+\Sigma d_i^{(M)}$ and ending at position number $1+k+\Sigma d_i^{(M)}$. And take $l$ to be any natural number such that $l\leq\Sigma d_{i+1}^{(M)}-\Sigma d_{i}^{(M)}-(k+1)$. 
\end{center}
\begin{center}
Or\\
For computing the possible values of the natural number $k$, proceed similarly as in sub case (I.a). And to determine the possible values of $l$, look at the entry in position number $\Sigma d_{i+s}^{(M)}+L+1$. If it is odd, then take $l=1$. If it is even, then look at the number of even entries in the string of binomial coefficients of $M$ that come in continuum starting from position number $1+\Sigma d_i^{(M)}$ (including the even entry at this position).  Let $p_j^{(M)}$ denote this number. Similarly, look at the number of even entries in the string of binomial coefficients of $M$ that come in continuum starting from position number $1+L+\Sigma d_{i+s}^{(M)}$ (including the even entry at this position).  Let $q_j^{(M)}$ denote this number. Tale $l$ to be any natural number such that $l\leq 1+min\{p_j^{(M)},q_j^{(M)}\}$.
\end{center}

\underline{\textbf{Case II:} When $\ ^{M}C_j$ is even.}\\
The position of $\ ^{M}C_j$ in the string of binomial coefficients of $M$ (starting from $\ ^{M}C_0$ to $\ ^{M}C_M$, that is, from left to right) is $(j+1)$-th.\\
\underline{Sub case II.a:} When $j+2=\Sigma d_i^{(M)}$ for some $i\in\{1,\ldots,2^m\}$ and there exists an odd natural number $u$ such that all the entries at the position numbers $\Sigma d_i^{(M)},1+\Sigma d_i^{(M)},\ldots, u+\Sigma d_i^{(M)}$ are odd and the entry at position number $u+1+\Sigma d_i^{(M)}$ is even.\\

I will now list down the various possibilities under this sub case:
\begin{center}
Either\\
Take $k$ to be any \textit{odd} natural number such that $k+1\leq u+1$. And take the string $\ ^{M}C_{j+1},\ldots,\ ^{M}C_{j+k+1}$ of length $(k+1)$ to be the continuous string of binomial coefficients of $M$ starting from position number $\Sigma d_i^{(M)}$ and ending at position number $k+\Sigma d_i^{(M)}$. Take $l$ to be any natural number such that $l\leq (u+1)-k$. 
\end{center}
\begin{center}
Or\\
Take the string $\ ^{M}C_{j+1},\ldots,\ ^{M}C_{j+k+1}$ of length $(k+1)$ to be the continuous string of binomial coefficients of $M$ starting from position number $j+2(=\Sigma d_i^{(M)})$ and ending at position number $\Sigma d_{i+s}^{(M)}$, where $s$ is any odd natural number such that $i+s\leq 2^m$. [Remark: Such a odd natural number $s$ may not exist always. If it does not exist, then this is an invalid possibility.] And to determine the possible values of $l$, look at the entry in position number $\Sigma d_{i+s}^{(M)}+1$. If it is even, then take $l=1$. If it is odd, then look at the number of odd entries in the string of binomial coefficients of $M$ that come in continuum starting from position number $\Sigma d_i^{(M)}$ (including the odd entry at this position).  Let $p_j^{(M)}$ denote this number. Similarly, look at the number of odd entries in the string of binomial coefficients of $M$ that come in continuum starting from position number $1+\Sigma d_{i+s}^{(M)}$ (including the odd entry at this position).  Let $q_j^{(M)}$ denote this number. Tale $l$ to be any natural number such that $l\leq 1+min\{p_j^{(M)},q_j^{(M)}\}$.
\end{center}
\noindent\underline{Sub case II.b:} When $j+2=\Sigma d_i^{(M)}$ for some $i\in\{1,\ldots,2^m\}$ and $\Sigma d_{i+1}^{(M)}-\Sigma d_{i}^{(M)}-1>0$.\\
Take the string $\ ^{M}C_{j+1},\ldots,\ ^{M}C_{j+k+1}$ of length $(k+1)$ to be the continuous string of binomial coefficients of $M$ starting from position number $j+2(=\Sigma d_i^{(M)})$ and ending at position number $\Sigma d_{i+s}^{(M)}$, where $s$ is any odd natural number such that $i+s\leq 2^m$. [Remark: Such a odd natural number $s$ may not exist always. If it does not exist, then this is an invalid possibility.] In this sub case, take $l=1$.\\
\noindent\underline{Sub case II.c:} When $j+2\neq\Sigma d_i^{(M)}$ for any $i\in\{1,\ldots,2^m\}$.
\begin{center}
Clearly there exists $i\in\{1,\ldots,2^m\}$ such that $j+2<\Sigma d_i^{(M)}$. Look at the smallest such $i$, call it $i_0$. Take the string $\ ^{M}C_{j+1},\ldots,\ ^{M}C_{j+k+1}$ of length $(k+1)$ to be the continuous string of binomial coefficients of $M$ starting from position number $j+2$ and ending at position number $\Sigma d_{i_{0}+s}^{(M)}$, where $s$ is any odd natural number such that $i_{0}+s\leq 2^m$. [Remark: Such a odd natural number $s$ may not exist always. If it does not exist, then this is an invalid possibility.] \\
And to determine the possible values of $l$, look at the entry in position number $1+\Sigma d_{i_{0}+s}^{(M)}$. If it is odd, then take $l=1$. If it is even, then look at the number of even entries in the string of binomial coefficients of $M$ that come in continuum starting from position number $j+2$ (including the even entry at this position).  Let $p_j^{(M)}$ denote this number. Similarly, look at the number of even entries in the string of binomial coefficients of $M$ that come in continuum starting from position number $1+\Sigma d_{i_{0}+s}^{(M)}$ (including the even entry at this position).  Let $q_j^{(M)}$ denote this number. Take $l$ to be any natural number such that $l\leq 1+min\{p_j^{(M)},q_j^{(M)}\}$.
\end{center}
There is \underline{\textit{another non-trivial possibility}} under this sub case, which is the following:
\begin{center}
Let $M=2^{N_r}+2^{N_{r-1}}+\cdots+2^{N_1}$ be the binary expansion of $M$ where $N_r>N_{r-1}>\cdots>N_1$ are non-negative integers. Let $N_{s_0}$ be the least element of the set $\{N_r, N_{r-1}, \cdots, N_1\}$ which is $\geq 2$. Let $\{N_{s_t}, N_{s_{t-1}}, \cdots, N_{s_1}\}$ be the subset of the set $\{N_r, N_{r-1}, \cdots, N_1\}$ defined by $\{N_b|1\leq b\leq r-1\ and\ N_{1+b}-N_{b}>1\}$. It may happen in some cases that $N_{s_0}=N_{s_1}$.

Given any $d\in\{1,\ldots,t\}$, let $\mathfrak{U}_d$ denote the set of all possible sums generated from elements of the set $\{2^{N_r},2^{N_{r-1}},\ldots,2^{N_{1+s_d}}\}$ (where each element of this set can appear atmost once in any such sum) which are strictly less than $2^{N_r}+2^{N_{r-1}}+\cdots+2^{N_{1+s_d}}$. Given any $d\in\{1,\ldots,t\}$ and any $x_{(d)}\in\mathfrak{U}_d$, suppose $\ ^{M}C_j$ is at the $(1+x_{(d)}+2^{1+N_{s_d}})$-th position in the string of binomial coefficients of $M$ \textit{counting from the end} (that is, from the right to left). In this situation, the following solutions are possible:\\

(i) Take $k+1$ to be equal to $2^{N_{s_d}}$. And to determine the possible values of $l$, look at the number of even entries in the string of binomial coefficients of $M$ that come in continuum after $\ ^{M}C_j$ (i.e., excluding $\ ^{M}C_j$ and to the right of it).  Let $p_j^{(M)}$ denote this number. Similarly, look at the number of even entries in the string of binomial coefficients of $M$ that come in continuum after $\ ^{M}C_{j+2^{(1+N_{s_d})}}$ (i.e., excluding $\ ^{M}C_{j+2^{(1+N_{s_d})}}$ and to the right of it).  Let $q_j^{(M)}$ denote this number. Tale $l$ to be any natural number such that $l\leq 1+2^{N_{s_d}}+min\{p_j^{(M)},q_j^{(M)}\}$.\\

(ii) If $d>1$ and there exists integer(s) $z$ such that $N_{s_d}>z>N_{s_{d-1}}$ and $z\in\{N_r, N_{r-1}, \cdots>N_1\}$, then for any such $z$, take $k+1$ to be equal to $2^{z}$. And to determine the possible values of $l$, look at the number of even entries in the string of binomial coefficients of $M$ that come in continuum after $\ ^{M}C_j$ (i.e., excluding $\ ^{M}C_j$ and to the right of it).  Let $p_j^{(M)}$ denote this number. Similarly, look at the number of even entries in the string of binomial coefficients of $M$ that come in continuum after $\ ^{M}C_{j+2^{(1+N_{s_d})}}$ (i.e., excluding $\ ^{M}C_{j+2^{(1+N_{s_d})}}$ and to the right of it).  Let $q_j^{(M)}$ denote this number. Tale $l$ to be any natural number such that $l\leq 1+2^{z}(2^{(N_{s_d}-z+1)}-1)+min\{p_j^{(M)},q_j^{(M)}\}$.\\

(iii) If $d=1$ and  $N_{s_1}>N_{s_{0}}$, then there always exists integer(s) $z$ such that $N_{s_1}>z\geq N_{s_{0}}$ and $z\in\{N_r, N_{r-1}, \cdots>N_1\}$. For any such $z$, do similarly as in (ii) above replacing $d$ there by $1$.
\end{center}
\subsection{Solution to Problem II}\label{ss.solutionprobII}
Problem II can be restated in the following way:---
\begin{center}
Given any two positive integers $\alpha$ and $\delta$ and any non-negative integer $j$, for what all values of positive integers $k,l$ and $q$ does the string
$$\ ^{M}C_j, \ ^{M}C_{j+1},\ldots,\ldots,\ ^{M}C_{j+l+k}$$
(where $M=\alpha+\delta$) of consecutive binomial coefficients of $M$ becomes a solution to Problem I as mentioned above \textit{and} the following properties are also satisfied (simultaneously):
\begin{itemize}
 \item $\ ^{M}C_{j+k+l+1}+\ ^{M}C_{j+k+l}+\cdots+\ ^{M}C_{j+l+1}$ and $\ ^{\alpha}C_{j+k+l+1}+\ ^{\alpha}C_{j+k+l}+\cdots+\ ^{\alpha}C_{j+l+1}$ should have the same parity.
 \item $\ ^{M}C_{j+k+l+2}+\ ^{M}C_{j+k+l+1}+\cdots+\ ^{M}C_{j+l+2}$ and $\ ^{\alpha}C_{j+k+l+2}+\ ^{\alpha}C_{j+k+l+1}+\cdots+\ ^{\alpha}C_{j+l+2}$ differ in parity if and only if $\ ^{\delta}C_1(\ ^{\alpha}C_{j+k+l+1}+\ ^{\alpha}C_{j+k+l}+\cdots+\ ^{\alpha}C_{j+l+1})$ is odd.
 \item $\ ^{M}C_{j+k+l+3}+\ ^{M}C_{j+k+l+2}+\cdots+\ ^{M}C_{j+l+3}$ and $\ ^{\alpha}C_{j+k+l+3}+\ ^{\alpha}C_{j+k+l+2}+\cdots+\ ^{\alpha}C_{j+l+3}$ differ in parity if and only if $\ ^{\delta}C_1(\ ^{\alpha}C_{j+k+l+2}+\ ^{\alpha}C_{j+k+l+1}+\cdots+\ ^{\alpha}C_{j+l+2})+\ ^{\delta}C_2(\ ^{\alpha}C_{j+k+l+1}+\ ^{\alpha}C_{j+k+l}+\cdots+\ ^{\alpha}C_{j+l+1})$ is odd.
 \item and so on till $\cdots\cdots$
 \item $\ ^{M}C_{j+k+l+q}+\ ^{M}C_{j+k+l+q-1}+\cdots+\ ^{M}C_{j+l+q}$ and $\ ^{\alpha}C_{j+k+l+q}+\ ^{\alpha}C_{j+k+l+q-1}+\cdots+\ ^{\alpha}C_{j+l+q}$ differ in parity if and only if $\Sigma_{p=1}^{q-1}\ ^{\delta}C_p(\ ^{\alpha}C_{j+k+l+q-p}+\ ^{\alpha}C_{j+k+l+q-1-p}+\cdots+\ ^{\alpha}C_{j+l+q-p})$ is odd.
\end{itemize} 
\end{center}
Given any positive integers $\alpha,\delta$ and non-negative integer $j$, we can proceed similarly as in subsection \ref{ss.solutionprobI} and find the possible values of the positive integers $k$ and $l$ for which the string
$$\ ^{M}C_j, \ ^{M}C_{j+1},\ldots,\ldots,\ ^{M}C_{j+l+k}$$
of consecutive binomial coefficients of $M$ becomes a solution to Problem I. But these values of positive integers $k$ and $l$ should also satisfy the additional $q$-many properties (about parities of sums of strings of length $(k+1)$ of the binomial coefficients of $M$ and $\alpha$) as mentioned above. 

For any two natural numbers $Z$ and $y$, let $\mathfrak{s}_{y}(Z)$ denote the sum of the first $y$-many (starting from $\ ^{Z}C_0$) binomial coefficients of $Z$. It is an easy exercise to check that for any such $Z$ and $y$, $\mathfrak{s}_{y}(Z)=\ ^{Z-1}C_{y-1}$. Since we are working over field of characteristic $2$, it is easy to see that $\ ^{M}C_{j+k+l+1}+\ ^{M}C_{j+k+l}+\cdots+\ ^{M}C_{j+l+1}$ equals $\mathfrak{s}_{j+k+l+1}(M)+\mathfrak{s}_{j+l}(M)$ which in turn equals $\ ^{M-1}C_{j+k+l}+\ ^{M-1}C_{j+l-1}$. Similarly one can say that $\ ^{\alpha}C_{j+k+l+1}+\ ^{\alpha}C_{j+k+l}+\cdots+\ ^{\alpha}C_{j+l+1}$ equals $\ ^{\alpha-1}C_{j+k+l}+\ ^{\alpha-1}C_{j+l-1}$. Therefore the condition that $\ ^{M}C_{j+k+l+1}+\ ^{M}C_{j+k+l}+\cdots+\ ^{M}C_{j+l+1}$ and $\ ^{\alpha}C_{j+k+l+1}+\ ^{\alpha}C_{j+k+l}+\cdots+\ ^{\alpha}C_{j+l+1}$ should have the same parity translates into the condition that
\begin{quote}
$\ ^{\alpha+\delta-1}C_{j+k+l}+\ ^{\alpha+\delta-1}C_{j+l-1}$ and $\ ^{\alpha-1}C_{j+k+l}+\ ^{\alpha-1}C_{j+l-1}$ should have the same parity. $\cdots\cdots\cdots$condition $(*)$
\end{quote}
The other $(q-1)$-many conditions (about parity of sums of strings of length $(k+1)$) imply that the positive integer $q$ should be such that for every $a\in\{2,\ldots,q\}$, the sums $\ ^{M}C_{j+k+l+a}+\ ^{M}C_{j+k+l+a-1}+\cdots+\ ^{M}C_{j+l+a}$ and $\ ^{\alpha}C_{j+k+l+a}+\ ^{\alpha}C_{j+k+l+a-1}+\cdots+\ ^{\alpha}C_{j+l+a}$ differ in parity if and only if the following condition holds:
\begin{quote}
If $p_1<\cdots<p_f$ is the collection of all elements of the set $\{p|p\in\{1,\ldots,a-1\}\ and\ \ ^{\alpha}C_{j+k+l+a-p}+\ ^{\alpha}C_{j+k+l+a-1-p}+\cdots+\ ^{\alpha}C_{j+l+a-p}\ is\ odd\}$, then \textit{exactly odd many} elements of the set $\{\ ^{\delta}C_{p_1},\ldots,\ ^{\delta}C_{p_f}\}$ should be odd. And this should hold true for every $a\in\{2,\ldots,q\}$.$\cdots\cdots\cdots$condition $(**)$
\end{quote}

We therefore need to determine when the value of the string sum $\ ^{M}C_{j+k+l+a}+\ ^{M}C_{j+k+l}+\cdots+\ ^{M}C_{j+l+a}$ becomes odd and when it is even as the integer $a$ ranges over the set $\{1,\ldots,q\}$. And similarly for the string sums $\ ^{\alpha}C_{j+k+l+a}+\ ^{\alpha}C_{j+k+l+a-1}+\cdots+\ ^{\alpha}C_{j+l+a}$ as the integer $a$ ranges over the set $\{1,\ldots,q\}$. But this problem is similar to solving problem I (or a problem equivalent to problem I where the column vector $b$ is replaced by a suitable vector which is either $(1,0,\ldots,0)^{t}$ or $(0,1,\ldots,1)^{t}$) for both the integers $M$ and $\alpha$, assuming (to begin with) that the string sums $\ ^{M}C_{j+k+l+1}+\ ^{M}C_{j+k+l}+\cdots+\ ^{M}C_{j+l+1}$ and $\ ^{\alpha}C_{j+k+l+1}+\ ^{\alpha}C_{j+k+l}+\cdots+\ ^{\alpha}C_{j+l+1}$ have the same parity. Recall that while solving problem I considering various sub cases, we got solutions for all possible values of the integer $l$. Speaking more precisely, we got all possible values of `upper bounds' of the natural number $l$. These upper bounds are nothing but the \textit{number of maximum possible strings (in continuum)} of length $(k+1)$ for which the parity of the string sum remains the same. That means, as soon as the value of $l$ exceeds this `upper bound', the string sum changes parity. Hence the positive integer $q$ should be such that the various upper bounds of the values of $l$ under consideration (for $\alpha$ and $M$ both) should tally with conditions $(*)$ and $(**)$ mentioned above. 
\section{Concluding remarks about rationality}\label{s.concluding-rem}
Looking at the rhythm of the parity of sums of continuous strings of binomial coefficients of any given positive integer $M$ (as described in subsections \ref{ss.solutionprobI} and \ref{ss.solutionprobII} above), one can depict that over fields of characteristic $2$, the Hilbert-Kunz multiplicity in this case of `disjoint-term trinomial hypersurfaces' will turn out to be rational. In fact, it will depend upon the positions of the odd and even entries in the entire string of binomial coefficients of $M$, a precise account of which is given in section \ref{s.structure-binomial-string}. I hope similar thing will hold true over fields of arbitrary positive characteristic $p$.

For any hypersurface defined by polynomials having `disjoint terms' in it (not just trinomial hypersurfaces), it is known that the corresponding Hilbert-Kunz multiplicity will be rational [see \cite{Mo-Teix2} and \cite{Mo-Teix1} for a reasoning]. The `phi´s´ in the papers of \cite{Mo-Teix2} and \cite{Mo-Teix1} attached to `disjoint-term trinomials´ (or more generally to any sum of monomials that are pairwise prime) are `$p$-fractals´. As a consequence, the Hilbert-Kunz series will be a rational function and the corresponding HK multiplicity will be in $\mathbb{Q}$. 

But for trinomial hypersurfaces NOT having `disjoint-terms' in it, the situation becomes more interesting. There we need to solve `similar' rank computation problems related to larger and more complicated systems of linear equations [in fact, there we can have infinitely many linear equations inside a single system]. But the solution to the rank computation problem mentioned in this article is like providing a basement for the work in the more general case. Due to the more complicated nature of the systems in that case, I suspect that the Hilbert-Kunz multiplicity can become irrational there.
\providecommand{\bysame}{\leavevmode\hbox
to3em{\hrulefill}\thinspace}
\providecommand{\MR}{\relax\ifhmode\unskip\space\fi MR }
\providecommand{\MRhref}[2]{%
  \href{http://www.ams.org/mathscinet-getitem?mr=#1}{#2}
} \providecommand{\href}[2]{#2}

\end{document}